\documentclass{amsart}

\usepackage{nicefrac}
\usepackage{mathrsfs}
\usepackage{xcolor}
\usepackage[colorlinks=true,linkcolor=colorref,citecolor=colorcita,urlcolor=colorweb]{hyperref}
\definecolor{colorcita}{RGB}{21,86,130}
\definecolor{colorref}{RGB}{5,10,177}
\definecolor{colorweb}{RGB}{177,6,38}

\newtheorem{theorem}{Theorem}[section]

\theoremstyle{definition}
\newtheorem{definition}[theorem]{Definition}

\theoremstyle{remark}

\numberwithin{equation}{section}

\newcommand{\N}{\mathbb{N}}

\begin{document}

\title{On sets of extreme functions for Fatou's theorem}

\author{Thiago R. Alves}
\address{Departamento de Matem\'{a}tica,
	Instituto de Ci\^{e}ncias Exatas,
	Universidade Federal do Amazonas,
	69.077-000 -- Manaus -- Brazil}
\email{alves@ufam.edu.br}
\thanks{The first author was supported in part by Coordenação de Aperfeiçoamento de Pessoal de Nível Superior - Brasil (CAPES) - Finance Code 001 and FAPEAM}

\author{Leonardo Brito}
\address{Departamento de Matem\'{a}tica,
	Instituto de Ci\^{e}ncias Exatas,
	Universidade Federal do Amazonas,
	69.077-000 -- Manaus -- Brazil}
\email{leocareiro2018@gmail.com}
\thanks{The second author was supported by FAPEAM}

\author{Daniel Carando}
\address{Departamento de Matem\'atica,
	Facultad de Cs. Exactas y Naturales,
	Universidad de Buenos Aires
	and IMAS-UBA-CONICET, Argentina.}
\email{dcarando@dm.uba.ar}
\thanks{The third author was supported by CONICET-PIP 11220130100329CO and  ANPCyT PICT 2018-04104.}

\subjclass[2020]{Primary 30H10, 32A35, 46B87; Secondary 46E25, 30H50}

\date{}

\dedicatory{}

\keywords{Fatou's theorem. Holomorphic functions. Hardy spaces. Algebrability. Spaceability.}

\begin{abstract}
Bounded holomorphic functions on the disk have radial limits in almost every direction, as follows from Fatou's theorem. Given a zero-measure set $E$ in the torus $\mathbb T$, we study the set of functions such that $\lim_{r \to 1^{-}} f(r \, w)$ fails to exist for every $w\in E$ (such functions were first constructed by Lusin). We show that the set of Lusin-type functions, for a fixed zero-measure set $E$, contain algebras of algebraic dimension $\mathfrak{c}$ (except for the zero function). When the set $E$ is countable, we show also in the several-variable case that the Lusin-type functions contain infinite dimensional Banach spaces and, moreover, contain plenty of $\mathfrak{c}$-dimensional algebras.
We also address the question for functions on infinitely many variables.
\end{abstract}

\maketitle

\section{Introduction and main results}

Throughout $\mathbb{D}$ will represent the open unit disc of the complex plane $\mathbb{C}$ and $\mathbb{T}$ will denote the torus in $\mathbb{C}$, i.e. $\mathbb{T}$ is the boundary of $\mathbb{D}$. A classical theorem by  Fatou \cite{Fatou} ensures  that for all bounded holomorphic function $f : \mathbb D \to \mathbb C$, the radial limits
$$\lim_{t \to 1} f(t \, w)$$
exists for almost every $w\in \mathbb T$. In other words, given a bounded holomorphic function $f : \mathbb D \to \mathbb C$, there exists a zero-measure $E\subset \mathbb T$ such that $\lim_{t \to 1} f(t \, w)$ exists for all $w\in \mathbb T\setminus E$. The limit in fact exists also \emph{non-tangentially}.
Lusin \cite{Lus} (see also \cite{Lus-Priv} and \cite[Chapter 2]{COLO66}) showed a kind of counterpart of Fatou's theorem: given a zero-measure set $E\subset \mathbb T$, there exists a  bounded holomorphic function $f : \mathbb D \to \mathbb C$ such that
\begin{equation}\label{eq-lusin}
\lim_{r \nearrow 1^-} f(r \, w) \text{  fails to exist for every  }w\in E.
\end{equation}
Our first result shows that the set of Lusin-type functions contains large algebras (except for the zero function). Recall that $H^\infty$  denotes the uniform algebra of bounded holomorphic functions on the unit disk.

\begin{theorem}\label{thm-lusin}
  Given a zero-measure $E\subset \mathbb T$, there exists a  subalgebra of $H^\infty$ of algebraic dimension $\mathfrak{c}$ every non-zero element of which satisfies \eqref{eq-lusin}.
\end{theorem}

Using the usual terminology, the previous theorem can be restated as follows: the set of ``Lusin-type" functions is \emph{strongly $\mathfrak{c}$-algebrable}. Let us recall the corresponding definition.

\begin{definition}
Let $\mathcal{A}$ be a commutative algebra and $\kappa$ be an infinite cardinal number. A subset $C \subset \mathcal{A}$ is said to be {\it strongly $\kappa$-algebrable} if $C \cup \{ 0 \}$ contains a subalgebra of $\mathcal A$ generated by an algebraically independent set of cardinality $\kappa$.
\end{definition}

The proof of Theorem \ref{thm-lusin}, which is a combination of Lusin's example with \cite[Theorem 7.5.1]{AronBerPelSeo}, is given in Section~\ref{sec-lusin}.
Before stating the next results, we present a main concept for this work. In the sequel, $B$ denotes  the unit ball of $\mathbb C^n$ with some norm and $S$ the corresponding unit sphere. Also, $H^\infty(B)$ stands for the algebra of bounded holomorphic functions on $B$.

\begin{definition}
	For $f\in H^\infty(B)$ and $x\in S$, we define the \emph{radial cluster set} of $f$ at $x$ as the set \begin{align}
	\nonumber rCl(f, x) = \{\mu \in \mathbb{C}\,: \mbox{ there exists a sequence } &(t_k)_{k=1}^{\infty} \subset (0,1) \mbox{ such that }\\ &  \lim_{k\to\infty}t_k=1 \ \ \textrm{and} \ \  \lim_{k\to\infty}f(t_k x)=\mu \}.
	\end{align}
\end{definition}
Note that the existence of radial limit at $x$ is equivalent to $rCl(f, x)$ being a singleton. Therefore, given a zero-measure set $E\subset \mathbb T$, Lusin-type functions are those having non-trivial cluster sets at each $w$ in  $E$.
Note also that we have the following expression for the radial cluster set:
	\begin{equation}\label{igualdade-cluster-linear}
		rCl(f,x)=\bigcap_{0<\delta<1}\overline{\{f(tx)\colon 1-\delta < t<1 \}}.
	\end{equation}
So, $rCl(f,x)$ can be seen as a decreasing intersection of connected compact sets and, therefore, is compact and connected (connectedness follows, for example, from  \cite[Exercise 11, Section 26]{Mun}). As a consequence, if $rCl(f,x)$ is not a singleton, it must have cardinality $\mathfrak{c}$, since it is a connected subset of $\mathbb C$.

On the other hand, since $rCl(f,x)$ is compact, we can consider it \emph{large} whenever it has non-empty interior. In this direction,  Matsushima \cite{Matsu} proved another kind of counterpart of Fatou's theorem in several variables: for any discrete subset $M$ of the Euclidean unit sphere $S^{n-1}$ of $\mathbb{C}^n$, $n \geq 1$, there is a bounded holomorphic function $f\in H^\infty(B_n)$ such that the radial cluster sets of $f$  at any point of $M$ contains  a fixed (closed) disc. In the previous notation, there is some disk $\mathbf D$ such that
$${\mathbf D}\subset \cap_{z \in M} \, rCl(f,z).$$
In particular, $f$ has large radial cluster set at every point of $M$.

The next results extends this result in several ways: we show large Banach subspaces and algebras in the set of functions which have big radial cluster set for every $x\in M$, for every countable set $M$ (not necessarily discrete). Moreover, the norm in $\mathbb C^n$ is arbitrary. The definitions needed for the next statements and their proofs are in be seen in Sections \ref{sec-2doteo-i} and \ref{sec-2doteo-ii}.

\begin{theorem}\label{thm1}
	Given a countable subset $M \subset S$,  let $\mathcal F\subset H^\infty(B)$ be the set of functions $f$ for which there exists a disc $\mathbf{D}$ such that
	\begin{eqnarray} \label{inclusion}
	\mathbf{D}\subset {\bigcap_{z \in M}} rCl(f,z).
	\end{eqnarray}
\begin{enumerate}
  \item[(I)] The set $\mathcal F$ contains (up to the zero function) an isomorphic copy of $\ell_1$. In particular, $\mathcal F$ is spaceable.
  \item[(II)] The set $\mathcal F$  is pointwise strongly $\mathfrak{c}$-algebrable and, in particular, strongly $\mathfrak{c}$-algebrable.
\end{enumerate}
\end{theorem}

\medskip
The  existence (almost everywhere) of radial limits for functions in $H^\infty$ leads to the well known isometric algebra isomorphism between $H^\infty$ and the Hardy space
\begin{equation*}
H^{\infty}(\mathbb{T}) = \{ g \in L_{\infty}(\mathbb{T}): \widehat{g}(n) = 0 \textit{ for  } n < 0  \}.
\end{equation*}
This isomorphism associates to each $f\in H^\infty(\mathbb{D})$ the function $f^*\in H^\infty(\mathbb T)$ given by
	$$f^*(e^{i\theta})=\lim_{t \nearrow 1^-} f(t \, e^{i\theta}).$$
Every function in $H^{\infty}(\mathbb{T})$ can be obtained in this way, taking $f$ as the convolution of $f^*$ with the Poisson kernel. For simplicity, we omit the $^*$ and write $f$ for both the holomorphic function on the disk and the function on the torus.

All this facts extend to Hardy spaces of several variables.
In \cite{Fatou-toru-infinito}, the authors studied to what extent Fatou's theorem holds for the Hardy spaces   on the
infinite-dimensional polytorus, in particular for the space
\begin{equation*}
H^{\infty}(\mathbb{T^{\infty}}) = \{ g \in L_{\infty}(\mathbb{T^{\infty}}): \widehat{g}(\alpha) \neq 0 \textit{ only if } \alpha_j \geq 0 \textit{ for every } j=1, 2, \ldots \}.
\end{equation*}
Functions in $H^{\infty}(\mathbb{T^{\infty}})$ correspond to holomorphic functions in $H^\infty(B_{c_0})$, via an isometric isomorphism of Banach algebras that preserve Fourier  and Taylor coefficients. Again, we write  $f$ for the both functions, i.e., $f\in  H^{\infty}(\mathbb{T^{\infty}})$ and  $f\in H^\infty(B_{c_0})$. In this context, the authors in \cite{Fatou-toru-infinito} study the problem of recovering the values of $f$ in $\mathbb{T^{\infty}}$ as radial limits. More precisely, given $f \in H^\infty(\mathbb T^\infty)$, one wonders if there exists a zero-measure (for the normalized Haar measure) set $E\subset \mathbb T^\infty$ such that for all $e^{i\theta} := (e^{i \theta_1}, e^{i \theta_2}, \ldots) \in \mathbb{T}^\infty \setminus E$ we have
	\begin{align} \label{EQ2}
	\lim_{n \to \infty} f(r_{1,n}e^{i\theta_1}, r_{2,n} e^{i\theta_2}, \ldots) = f(e^{i\theta})
	\end{align}
whenever $(r_{k,n})_{k=1}^\infty \in c_0 \cap [0,1]^\infty$ for every $n$ with $\lim\limits_{n \to \infty} r_{k,n} = 1$ for every $k$.

The authors in \cite{Fatou-toru-infinito} show a counter example showing that answers the previous question by the negative. Moreover, in their example  \eqref{EQ2} fails almost everywhere. Our last result considers  the set of such counter examples and reads as follows.

\begin{theorem} \label{thm1.5}
	There is a zero-measure subset $F$ of $\mathbb{T}^\infty$ such that the set of functions in $H^\infty(\mathbb{T}^\infty)$ for which (\ref{EQ2}) fails for every $e^{i\theta} \in \mathbb{T}^\infty \setminus F$ is strongly $\mathfrak{c}$-algebrable.
\end{theorem}

\section{Proof of Theorems \ref{thm-lusin} and \ref{thm1.5}} \label{sec-lusin}

The proof of Theorem \ref{thm-lusin} is rather simple once we have a function $f$ satisfying \eqref{eq-lusin}. By  \cite[Theorem 7.5.1]{AronBerPelSeo}, it is enough to show that, for every $\varphi$ of exponential type, the composition $\varphi\circ f$ also satisfies \eqref{eq-lusin}. We say that $\varphi$ is of exponential type if it can be written as $\varphi(z)=\sum_{j=1}^{N}a_{j}\exp(b_{j}z)\in\mathcal{E}$, with $N\in\mathbb{N}$, $a_{1},\ldots,a_{N},b_{1},\ldots,b_{N}\in\mathbb{C}\setminus\{0\}$,  and $b_{1},\ldots,b_{N}$ all different. Just of completeness, let us mention that if $B\subset \mathbb{R}$ is a $\mathbb Q$-linear independent set of cardinality $\mathfrak{c}$, then  $B_f=\{e^{bf}\colon b\in B\}$ can be seen to be an algebraically independent subset of $H^\infty$. Also, every function in the algebra generated by $B_f$ is of the form $\varphi\circ f$ with $\varphi$ of exponential type. Therefore, if we show that $\varphi\circ f$ satisfies \eqref{eq-lusin}  for every such $\varphi$, we are done.

Let $w\in E$ be fixed. For $\mu\in rCl(f,w)$ we can take $(t_k)_{k\in \mathbb N}\subset(0,1)$ such that $\lim\limits_{k\to\infty}t_k=1$ and $\lim\limits_{k\to\infty}f(t_k w)=\mu$. Therefore, since $\varphi$ is continuous in the whole plane,
\begin{align}\label{tipo-exponecial-cluster-liner}
\lim_{k\to\infty}(\varphi\circ f)(t_k w)= \varphi (\mu).
\end{align}
This means that $\varphi(rCl(f,w))$ is contained in $rCl(\varphi\circ f, w)$. Now, $\varphi$ is a non-constant entire function (this can be seen in the proof of \cite[Theorem 7.5.1]{AronBerPelSeo}). Since $rCl(f,w)$ is a compact set of cardinality $\mathfrak c$, $\varphi$ cannot be constant on $\mathfrak c$ by the identity principle. Therefore,  $\varphi(rCl(f,w))$ has more than one point and, by the inclusion above, so has $rCl(\varphi\circ f, w)$. Since this holds for any $w\in E$, we conclude that $\varphi\circ f$ satisfies \eqref{eq-lusin} and this proves Theorem \ref{thm-lusin}.

\medskip
 Let us now prove Theorem \ref{thm1.5}, which follow similar lines.
From \cite[Theorem 4(i)]{Fatou-toru-infinito}, there exist a function $f \in H^\infty(\mathbb{T}^\infty)$, a zero-measure set $F_1$ of $\mathbb{T}^\infty$, and a sequence ${\bf r}_n := (r_{k,n})_{k=1}^\infty$ in $c_0 \cap [0,1]^\infty$ with $\lim\limits_{k\to\infty}r_{k,n}=1$, so that $f(e^{i\theta}) \not= 0$ for every $e^{i\theta} \in \mathbb{T}^\infty$ and
	\begin{align} \label{EQ-2-torus}
	\lim_{n \to \infty} f(r_{1,n} e^{i \theta_1}, r_{2,n} e^{i \theta_2}, \ldots) = 0
	\end{align}
	for every $e^{i \theta} = (e^{i \theta_1}, e^{i \theta_2}, \ldots) \in \mathbb{T}^\infty \setminus F_1$.
On the other hand, by \cite[Theorem 1(iii)]{Fatou-toru-infinito}, there exists another zero-measure set $F_2$ of $\mathbb{T}^\infty$ such that
	\begin{align} \label{EQ-1-torus}
	\lim_{r \to 1^-}f(r e^{i\theta_1}, r^2 e^{i \theta_2}, r^3 e^{i \theta_3}, \ldots) = f(e^{i \theta})
	\end{align}
	whenever $e^{i \theta} := (e^{i \theta_1}, e^{i \theta_2}, \ldots) \in \mathbb{T}^\infty \setminus F_2$.
Setting $F=F_1\cup F_2$, this means that for each $e^{i \theta} \in \mathbb{T}^\infty \setminus F$ the set
	\begin{align*}
		\nonumber   Cl_\infty(f, e^{i\theta}) := \{\mu \in \mathbb{C}\,:& \mbox{ there is a sequence } (\mbox{\bf r}_n)_{n=1}^\infty\subset c_0 \cap [0,1]^\infty \mbox{ such that }\\ &  \lim\limits_{k\to\infty}r_{k,n}=1 \ \ \textrm{and} \ \  \lim\limits_{n\to\infty}f(r_{1,n}e^{i\theta_1},r_{2,n}e^{i\theta_2},\ldots) = \mu \}
	\end{align*}
	 has at least two elements. In fact, it must have the cardinality of the continuum, since we can again see that $Cl_\infty(f, e^{i\theta}) $ is a connected subset of the complex plane. Indeed, we consider $K := \Pi_{k=1}^\infty [0,1]e^{i \theta_k}$ endowed with the product topology and take  $(U_\lambda)_{\lambda\in \Lambda}$ a decreasing net that form a basis of connected neighbourhoods of $e^{i \theta}$ in $K$ (we can take $\lambda=(A, \varepsilon)$ for $A\subset \mathbb N$ and $\varepsilon>0$ with the natural order and define $U_\lambda$ in the natural way). With this, we get
	 \begin{align*}
	 	Cl_\infty(f,e^{i\theta}) = \bigcap_{\lambda \in \Lambda} \overline{f(U_\lambda \cap B_{c_0})}.
	 \end{align*}
That is, $Cl_\infty(f,e^{i\theta})$ may be written as a decreasing intersection of connected compact subsets. Hence, it follows from \cite[Exercise 11, Section 26]{Mun} that the set $Cl_\infty(f,e^{i\theta})$ must be compact and connected as well.

Hence, we fix the function $f \in H^\infty(B_{c_0})$  above, for which  $ Cl_\infty(f,e^{i \theta})$ has cardinality $\mathfrak c$ for each $e^{i \theta} \in \mathbb{T}^\infty \setminus E$. Now,  we may conclude the proof of Theorem \ref{thm1.5} proceeding as in the proof of Theorem $\ref{thm-lusin}$, making use of \cite[Theorem 7.5.1]{AronBerPelSeo}.

\section{Proof of part (I) of Theorem \ref{thm1}} \label{sec-2doteo-i}

Recall that a sequence $(x_k)_k$ in $B$ is an \emph{interpolating sequence} for $H^\infty(B)$ if for any sequence $(\lambda_k)_k \in \ell_\infty$, there is ${f} \in H^\infty(B)$ such that ${f}(x_k) = \lambda_k$. In that case, there exists a constant $C\ge 1$ (the \emph{interpolating constant of} $(x_k)_k$) such that we can always get $f$ satisfying also $\|f\|_\infty \le C \|(\lambda_k)_k\|_\infty$.

In  \cite{GalMir} the authors extend the so-called Hayman-Newman condition to arbitrary Banach spaces. In our context, their result reads as follows: if there is some $0 < c < 1$ such that \begin{equation}\label{eq-H-N}
                                                  (1-\|x_{k+1}\|)/(1-\|x_k\|) < c \text{ for all }k,
                                                \end{equation} then
	 $(x_k)_k$ is an interpolating sequence for ${H}^\infty(B)$.

From now on, we fix an enumeration $\{s_q : q \in \mathbb{N}\}$ of the rational numbers $\mathbb{Q}$  and denote by  $(\mathfrak{p}_{k})_{k=1}^{\infty}$ the sequence of prime numbers. Let us start with the proof of the spaceability part in Theorem \ref{thm1}.

Given a countable subset $M\subset S$, we start by defining an interpolating sequence such that for each $z\in M$ there is a subsequence converging radially to $z$. To this end, we write $M=\{z_\ell:\ell\in \N\}$ and we take a partition $\{\Theta_\ell : \ell \in \mathbb N\}$ of the natural numbers into (disjoint) infinite subsets. Now, for $m\in \Theta_\ell$ we set $x_{m}=t_m^\ell z_\ell$ for some $t_m^\ell \in (0,1)$. If we take $t_m^\ell$ converging to 1 fast enough (as $m\to\infty$), we can ensure that $(x_{m})_m$ satisfies condition \eqref{eq-H-N}. That is, we can take $(x_{m})_m$ to be an interpolating sequence.
Note also that $(x_m)_{m\in \Theta_\ell}$  converges to $z_\ell$ for  $m\to\infty, \ m\in \Theta_\ell$.

Let us, for each $\ell\in \N$, take a (disjoint) partition $\{\Theta_{\ell,n}:n\in\mathbb{N}\}$ of $\Theta_\ell$. We write $$\Theta_{\ell,n} := \{m_{\ell,n,1}<m_{\ell,n,2}< \cdots \}.$$ Let  $(a_k)_{k=1}^{\infty} \subset \mathbb{D}$ be a sequence which accumulates at every point of $\overline{\mathbb D}$. Since $(x_{m})_m$ is an interpolating sequence, for each $j\in \N$ we can choose a function $f_j\in {H}^\infty(B)$ such that
$${f}_j(x_{m_{\ell,n,q}}) = \mathfrak{p}_j^{is_{q}} \cdot a_{n},\quad\text{ for each } \ell,n,q\in \N$$ and such that  $\displaystyle\sup_{j\in\mathbb{N}}{\Vert f_j\Vert}_\infty = C<\infty$.

We define $T:\ell_1\longrightarrow {H}^{\infty}(B)$ by $$T(\lambda)=\sum_{j=1}^\infty \lambda_j f_j, \quad \lambda=(\lambda_j)_{j}.$$ It is clear that $T$ is a bounded operator with $\|T\|\le C$. To see that it is an isomorphism with its image, we fix $\varepsilon>0$  and take $M_0$ such that $\sum_{j>M_0}|\lambda_j|<\varepsilon/2$. This gives
\begin{eqnarray}
  \Big\Vert \sum_{j=1}^\infty \lambda_j f_j\Big\Vert_\infty &\ge& \sup_{\ell,n,q\in \N} \Big| \sum_{j=1}^\infty \lambda_j {f}_j(x_{m_{\ell,n,q}}) \Big| = \sup_{n,q\in \N} \Big| \sum_{j=1}^\infty \lambda_j \mathfrak{p}_j^{is_{q}} \cdot a_{n} \Big| \nonumber \\
   &=&  \sup_{q\in \N} \Big| \sum_{j=1}^\infty \lambda_j \mathfrak{p}_j^{is_{q}}  \Big| \ge \sup_{q\in \N} \Big| \sum_{j=1}^{M_0} \lambda_j \mathfrak{p}_j^{is_{q}} \Big|-\frac{\varepsilon}{2}. \label{eq-bdedbelow}
\end{eqnarray}
 By the density properties of  Kronecker flow in $\mathbb T^{M_0}$ (see, for example, \cite[Proposition 3.4]{DefGarMaeSev19}), we have
 $$ \sup_{q\in \N} \Big| \sum_{j=1}^{M_0} \lambda_j \mathfrak{p}_j^{is_{q}}  \Big| = \sum_{j=1}^{M_0} |\lambda_j|\ge\sum_{j=1}^{\infty} |\lambda_j| -\frac\varepsilon 2 .$$
This, together with  \eqref{eq-bdedbelow}, shows that $T$ is bounded below and then it is an isomorphism with its image. Now we show that for $\lambda\ne 0$ and $f=T(\lambda)$ we have
\begin{equation}\label{eq-spaceab}
  \|\lambda\| \overline{\mathbb D}\subset {\bigcap_{\ell\in\N}} rCl(f,z_\ell).
\end{equation}
Again, we fix $\varepsilon>0$ and now we take $M_0$ such that $\sum_{j>M_0}|\lambda_j|<\varepsilon/4$. By the density of Kronecker flow, we can choose infinitely many values of $q$ such that
$$  \Big| \sum_{j=1}^{M_0} \lambda_j \mathfrak{p}_j^{is_{q}}  - \sum_{j=1}^{M_0} |\lambda_j| \Big| <\frac\varepsilon 4 .$$
Finally, given $\mu\in \|\lambda\|\overline{\mathbb D}$, there are infinitely many $n$ such that $\big| \|\lambda\| a_n-\mu \big|<\varepsilon/4$. For such $q$ and $n$ and for any $\ell\in\N$ we have
\begin{eqnarray}
\nonumber\Big| \sum_{j=1}^\infty \lambda_j {f}_j(x_{m_{\ell,n,q}}) - \mu\Big|&\le & \Big| \sum_{j=1}^{M_0} \lambda_j {f}_j(x_{m_{\ell,n,q}}) - \mu\Big| +\frac{\varepsilon}{4} = \Big| \sum_{j=1}^{M_0} \lambda_j \mathfrak{p}_j^{is_{q}} \cdot a_{n}  - \mu\Big| +\frac{\varepsilon}{2} \\
\nonumber   &\le&    \Big| \sum_{j=1}^{M_0} |\lambda_j | \cdot a_{n}  - \mu\Big| +\frac{\varepsilon}{4} +\frac{\varepsilon}{4} \\
   &\le& \Big| \sum_{j=1}^{\infty} |\lambda_j | \cdot a_{n}  - \mu\Big| +\frac{\varepsilon}{4}+\frac{\varepsilon}{4} +\frac{\varepsilon}{4} <\varepsilon. \label{eq-mu}
\end{eqnarray}
Recall that, by construction,  $x_{m_{\ell,n,q}} = t^\ell_{m_{\ell,n,q}} z_\ell$ and $t^\ell_{m_{\ell,n,q}} \to 1 $ as $n,q\to\infty$. This and \eqref{eq-mu} give, for any $\delta>0$, some $t=t^\ell_{m_{\ell,n,q}}>1-\delta$ such that $|f(tz_\ell)-\mu|<\varepsilon$. Thus, $\mu$ belongs to $rCl(f, z_\ell)$ for every $\ell\in\N$. This shows
 that every non-zero function $f$ in the image of $T$ satisfies \eqref{eq-spaceab} for $\lambda=T^{-1}(f)$, as desired.

\section{Proof of part (II) of Theorem \ref{thm1}} \label{sec-2doteo-ii}

For $\mathcal{A}$ a commutative algebra and $C$  a non-empty subset of $\mathcal{A}$, we say that $C$ is {\it pointwise strongly $\frak{c}$-algebrable} if for every $f \in C$ there is a subalgebra $\mathcal{B}$ of $\mathcal{A}$ which is generated by an algebraically independent set $G$ with $\#G = \frak{c}$ and $f \in \mathcal{B} \subset C \cup \{0\}$. It is worth noting that both the terminology and the notion of pointwise strongly algebrability is inspired by previous works concerning a similar notion for lineability/spaceability (see \cite{FavPelTom, PelRap}).

Next we prove the pointwise strong algebrability part of Theorem \ref{thm1}. Choose $f \in H^\infty(B)$ satisfying
\begin{eqnarray} \label{inclusion2}
	\mathbf{D} \subset {\bigcap_{z \in M}} rCl(f,z)
	\end{eqnarray}
 for some disk $\mathbf{D}$. We may define an interpolating sequence $(x_m)$ for $H^\infty(B)$ such that for each $z \in M$ there is a subsequence $(x_{m_k})_k$ of $(x_m)$ converging radially to $z$ and for which $(f(x_{m_k}))_k$ accumulates at every point of $\mathbf{D}$. Indeed, we write $M = \{z_\ell : \ell \in \mathbb N\} \subset S$ and proceed similarly as we did in the proof of the spaceability part of Theorem \ref{thm1}.  We take sequences of the form $s_m^\ell \, z_\ell$, with  $s^\ell_m$ converging to $1$ fast enough (as $m \to \infty$) so that together they form an interpolating sequence, and so that  ${(f(s_m^\ell \, z_\ell))}_m$ accumulates at every point of $\mathbf{D}$, which is possible because $f$ satisfies (\ref{inclusion2}).

Let us now fix some notations. First, let $(a_p) \subset \mathbf{D}$ be a sequence that accumulates at every point of  $\mathbf{D}$. For each $(\ell, p) \in \mathbb{N}^2$, let $\Theta^{\ell,p} := \{m_1^{\ell,p} < m_2^{\ell,p} < m_3^{\ell,p} < \cdots\} \subset \mathbb N$ be a set for which the sequence $(x_{m_j^{\ell,p}})_j$ converges radially to $z_\ell$ and $(f(x_{m_j^{\ell,p}}))_j$ converges to $a_p$. We may and do take $\{\Theta^{\ell,p}\}_{\ell,p}$ as a (pairwise disjoint) partition of the set of all natural numbers. Also, for each $(\ell, p) \in \mathbb{N}^2$, we let $\{\Theta_q^{\ell,p}\}_{q \in \mathbb{N}}$ be a (pairwise disjoint) partition of infinite subsets of $\Theta^{\ell,p}$; say, $\Theta_q^{\ell,p} := \{m_{q,1}^{\ell,p} < m_{q,2}^{\ell,p} < m_{q,3}^{\ell,p} < \cdots\}$.

Fix a bijection $u : \mathbb{Q}_+^* \to \mathbb N$. Since the collection of sets $\Theta_q^{\ell,p}s$ are pairwise disjoint and cover the set of all natural numbers, we get a well-defined bijection
$w : \mathbb{N}^3 \times \mathbb{Q}_+^* \to \mathbb N $ given by $w(j,\ell,p,r) :=   m_{u(r),j}^{\ell,p}$. To simplify notation, we set $y_{r,j}^{\ell,p} := x_{w(j,\ell,p,r)}$ for $r \in \mathbb{Q}_+^*$ and $j,\ell,p \in \mathbb{N}$. Note that $(y_{r,j}^{\ell,p})_j$ is a subsequence of $(x_{m_j^{\ell,p}})_j$ for every $\ell, p \in \mathbb{N}$ and every $r \in \mathbb{Q}_+^*$. In particular, $(y_{r,j}^{\ell,p})_j$ converges radially to $z_\ell$ and $f(y_{r,j}^{\ell,p}) \to a_p$ as $j \to \infty$.

Now, for every $k\in\mathbb{N}$ and every $\xi \in (0,1)$, we take $k-1=t^{k}_{\xi,0}<t^{k}_{\xi,1} < \cdots<t^{k}_{\xi,N^{k}_{\xi}}=k$ a partition of the interval $[k-1,k]$ such that $\vert t^{k}_{\xi,i}-t^{k}_{\xi,i-1}\vert < \xi^{k}$ and $N^{k}_{\xi}\in\mathbb{N}$. We set $I^{k}_{\xi,i}:=[t^{k}_{\xi,i-1},t^{k}_{\xi,i})$ for $k\in\mathbb{N}$, $\xi\in(0,1)$ and $i=1,2,\ldots,N^{k}_{\xi}$. Since $(x_m)$ is an interpolating sequence for $H^{\infty}(B)$, for each $\xi \in (0,1)$ there is $g_\xi \in H^{\infty}(B)$ so that
\begin{eqnarray}\label{Func.g_xi}
	g_{\xi}(y_{r,j}^{\ell,p}) := \begin{cases}
		f(x_{m_{j}^{\ell,p}}) & \text{ if } r \in I^{k}_{\xi,i} \mbox{ and } i \equiv 2\ (\mathrm{mod}\ 3);
		\\
		1 & \text{ if } r\in I^{k}_{\xi,i} \mbox{ and } i\equiv 1\ (\mathrm{mod}\ 3);
		\\
		0 & \textrm{otherwise}.
	\end{cases}
\end{eqnarray}

From now on, our proof follows similar lines to that of \cite[Lemma 2.6]{AlvCar} and, thus, we omit some details. First, we set $C_0^{j,\ell,p} := 0$, $C_1^{j,\ell,p} := 1$ and $C_2^{j, \ell, p} := f(x_{m_j^{\ell,p}})$ for every $j, \ell, p \in \mathbb N$. It follows from (\ref{Func.g_xi}) that, for given $N \in \mathbb{N}$, $\xi_1 < \xi_2 < \cdots < \xi_N$ in $(0,1)$ and $(i_1, \ldots, i_N) \in \{0,1,2\}^N$, we can find $k_0 \in \mathbb{N}$ and $r_0 \in [k_0 - 1, k_0) \cap \mathbb{Q}$ satisfying
\begin{eqnarray} \label{eq-proof}
	(g_{\xi_1}(y_{r_0,j}^{\ell,p}), \ldots, g_{\xi_N}(y_{r_0,j}^{\ell,p})) = (C_{i_1}^{j,\ell,p}, \ldots, C_{i_N}^{j,\ell,p})
\end{eqnarray}
for every $j, \ell, p \in \mathbb{N}$.
	Set $g_1 := f$. To complete the proof we  show that the set $\{g_\xi : \xi \in (0,1]\}$ is an algebraically independent set which generates a sub-algebra $\mathcal{B}$ lying in $\mathcal{F} \cup \{0\}$. Similarly as in the proof of \cite[Lemma 2.6]{AlvCar}, this follows from the assertion below:
		\begin{enumerate}
		\item[$(*)$] For every $N \in \mathbb N$, each choice of numbers $\xi_1 < \cdots < \xi_N$ in $(0,1]$ and every non-constant polynomial $Q \in \mathbb C[X_1, \ldots, X_N] \setminus \{0\}$, there is a disc $\mathbf{D'}\subset rCl(Q(g_{\xi_1}, \ldots, g_{\xi_N}), z_\ell)$ for every $\ell$.
	\end{enumerate}
Indeed, $(*)$ implies that $\mathcal{B}$ is contained in $\mathcal{F} \cup \{0\}$ and, also, that no algebraic combination of $g_\xi$'s may be zero.

Let us prove Assertion $(*)$. For $Q$ as in the statement, let $P$ be the polynomial on $N$ complex variables which comes from $Q$ after turning all non-trivial powers of $X_1,\dots, X_{N-1}$ to 1. For example, if $Q=X_1^2 X_2+X_1^3 X_2^2+   X_3^3X_N^5$, then $P=2 X_1 X_2+   X_3X_N^5$. The important feature of $P$ is that $P(X)=Q(X)$ whenever $X_j=  0$ or $X_j=1$ for $j=1, \dots, N-1$.

Given $A\subset\{1,\dots,N-1\}$ we write $X_A=\prod_{j\in A} X_j$. With this notation, we can write
\begin{equation}\label{eq-defdeP}
P(X)= \sum_{k=0}^{d} \Big(\sum_{A\subset \{1,\dots,N-1\}} C_{A,k} X_A \Big) X_N^k,
\end{equation}
where $d$ is the maximal power of $X_N$ in $Q$. We now take $A_0\subset \{1,\dots,N-1\}$ such that $C_{A_0,d}\ne 0$ and that $\# A_0\le \#A$ for all $A\subset \{1,\dots,N-1\}$ with $C_{A,d}\ne 0$. Such $A_0$ clearly exists and has the following property: if $C_{A,d}\ne 0$, then $A\setminus A_0$ is non-empty. Equivalently, if  $X_j=1$ for $j\in A_0$ and $X_j=0$ otherwise, then $X_A=0$ for all $A$ with $C_{A,d}\ne 0$, $A\ne A_0$. And, of course, $X_{A_0}=1$.

Now, for every $j, \ell, p \in \mathbb{N}$ we take $i_k = 1$ if $k\in A_0$ and $i_k = 0$ if $k\in \{1,\dots, N-1\}\setminus A_0$. By using (\ref{eq-proof}) we can find $r_0 \in \mathbb{Q}_+$ such that   $g_{\xi_k}(y_{r_0,j}^{\ell, p})=1$ for $k\in A_0$ and $g_{\xi_k}(y_{r_0,j}^{\ell, p})=0$ otherwise (for every $j, \ell, p \in \mathbb{N}$, $k=1,\dots,N-1$). Note that for $X=(g_{\xi_1}(y_{r_0,j}^{\ell, p}), \ldots,g_{\xi_N}(y_{r_0,j}^{\ell,p}))$ and each $A\subset \{1,\dots,N-1\}$ we have $X_A=0$ or $1$ and this value is independent of $j, \ell, p \in \mathbb{N}$.
As a consequence, using the expression \eqref{eq-defdeP} we have
\begin{multline*}
   Q(g_{\xi_1}(y_{r_0,j}^{\ell, p}), \ldots,g_{\xi_N}(y_{r_0,j}^{\ell,p}))   = P(g_{\xi_1}(y_{r_0,j}^{\ell, p}), \ldots,g_{\xi_N}(y_{r_0,j}^{\ell,p})) \\
  =  C_{A_0,d} \,\,  g_{\xi_N}(y_{r_0,j}^{\ell,p})^d  + \sum_{k=0}^{d-1} \Big(\sum_{A\subset \{1,\dots,N-1\}} C_{A,k} X_A \Big) \,g_{\xi_N}(y_{r_0,j}^{\ell,p})^d \\
  =  C_{A_0,d} \,\,  g_{\xi_N}(y_{r_0,j}^{\ell,p})^d + \text{ lower powers of }g_{\xi_N}(y_{r_0,j}^{\ell,p}) \\
  \hspace{-5cm} =  \widetilde q(g_{\xi_N}(y_{r_0,j}^{\ell,p})),\hspace{6cm} \,\,
\end{multline*}
for $\tilde q$ some one variable polynomial of degree $d\ge 1$ depending only on $Q$. By the open mapping theorem, $\widetilde{q}(\mathbf D)$ contains a disk $\mathbf{D'}$. Since for every $\ell,p \in \mathbb{N}$ we have that $(x_{m_j^{\ell,p}})_j$ converges radially to $z_\ell$, $f(x_{m_j^{\ell,p}}) \to a_p$ as $j \to \infty$, and $(y_{r_0,j}^{\ell,p})_j$ is a subsequence of $(x_{m_j^{\ell,p}})_j$, then the disc $\mathbf{D'}$ must be contained in $rCl(Q(g_{\xi_1}, \ldots, g_{\xi_N}), z_\ell)$ for each $\ell$. This proves Assertion $(*)$.

\bibliographystyle{amsplain}

\end{document}